\documentstyle[12pt]{article}

\font\tensym=msbm10
\font\sevensym=msbm7
\font\fivesym=msbm5

\font\tengoth=eufm10
\font\sevengoth=eufm7
\font\fivegoth=eufm5

\newfam\symfam
\textfont\symfam=\tensym
\scriptfont\symfam=\sevensym
\scriptscriptfont\symfam=\fivesym
\def\sym{\fam\symfam\tensym}

\newfam\gothfam
\textfont\gothfam=\tengoth
\scriptfont\gothfam=\sevengoth
\scriptscriptfont\gothfam=\fivegoth
\def\goth{\fam\gothfam\tengoth}

\newtheorem{example}{Example}[section]

\newtheorem{theorem}[example]{Theorem}
\newtheorem{corollary}[example]{Corollary}

\newtheorem{proposition}[example]{Proposition}

\def\boxit#1#2{\setbox1=\hbox{\kern#1{#2}\kern#1}%
\dimen1=\ht1 \advance\dimen1 by #1 \dimen2=\dp1 \advance\dimen2 by #1
\setbox1=\hbox{\vrule height\dimen1 depth\dimen2\box1\vrule}%
\setbox1=\vbox{\hrule\box1\hrule}%
\advance\dimen1 by .4pt \ht1=\dimen1
\advance\dimen2 by .4pt \dp1=\dimen2 \box1\relax}

\def\Proof{\medskip\noindent {\it Proof --- \ }}
\def\cqfd{\hfill $\Box$ \bigskip}
\def\adots{\mathinner{\mkern2mu\raise1pt\hbox{.}
\mkern3mu\raise4pt\hbox{.}\mkern1mu\raise7pt\hbox{.}}}
\def\<{\langle}
\def\>{\rangle}

\def\ie{{\it i.e.\,}}

\def\carrel#1#2#3{{\vcenter{\hrule height .#2pt
                       \hbox{\vrule height #1pt width .#2pt \kern #3pt
                             \vrule width .#2pt}
                       \hrule height .#2pt}}}
\def\carre{\mathchoice\carrel87{8.7}\carrel87{8.7}
                      \carrel55{5.5}\carrel45{4.5} \, }

\def\S{\hbox{\goth S}}
\def\Sym{{\bf Sym}}

\def\N{{\sym N}}
\def\C{{\sym C}}
\def\Z{{\sym Z}}

\def\u{{\bf u}}
\def\x{{\bf x}}

\title{\bf Flag Varieties and the Yang-Baxter Equation }

\author{Alain {\sc Lascoux}\thanks{Institut Blaise pascal, LITP,
Universit\'e Paris 7, 2 Place Jussieu, 75251 Paris Cedex 05, France}, \
Bernard {\sc Leclerc}\thanks{Universit\'e de Caen,  D\'epartement de
 Math\'ematiques, Esplanade de la Paix, BP 5186, 14032 Caen Cedex, France}\
\ and \  
Jean-Yves {\sc Thibon}\thanks{
Institut Gaspard Monge,
Universit\'e de Marne-la-Vall\'ee,
2, rue de la Butte-Verte,
93166 Noisy-le-Grand Cedex, France} 
}

\date{}

\begin{document}

\newfont{\rml}{cmu10}
\newfont{\rmL}{cmr12}
\maketitle
\date{}

\begin{abstract}
We investigate certain bases of  Hecke algebras defined
by means of the Yang-Baxter equation, which we call Yang-Baxter bases. 
These bases are essentially
self-adjoint with respect to a canonical bilinear form. 
In the case of the degenerate Hecke algebra,
we identify the coefficients in the expansion of the Yang-Baxter basis
on the usual basis of the algebra 
with specializations of double Schubert polynomials.
We also describe the
expansions associated to other specializations of the generic Hecke algebra.
\end{abstract}

\section{Introduction}

Yang's original motivation for introducing the Yang-Baxter equation \cite{Y}
was the $n$-body problem on a circle with Hamiltonian
$$
H(\x)= -\sum_{i=1}^n{\partial^2\over \partial x_i^2} + 2c\sum_{i<j}\delta(x_i-x_j) \ ,
$$
where $\delta$ is the Dirac distribution.
The problem was to solve the Schr\"odinger equation
$$
H(\x)\psi(\x)=E \psi(\x)
$$
with periodic boundary conditions.
Using a variant of the Bethe Ansatz, Yang looked for solutions of the form

$$\psi(\x) = \sum_{\tau\in\S_n}  \theta(\x^\tau) \sum_{\mu\in\S_n} A_\mu(\tau) 
e^{{\rm i}\u^\mu \cdot \x^\tau}$$
where $\S_n$ is the symmetric group,
$\u =(u_1,\ldots,u_n)\in\C^n$ is the vector of spectral parameters,
$\u^\mu=(u_{\mu(1)},\ldots,u_{\mu(n)})$,
and $\theta$ is the characteristic function of the domain
$x_1<x_2<\ldots <x_n$. 
The unknown coefficients $A_\mu(\tau)$ form an $n!\times n!$ matrix, and
it is convenient to regard each $A_\mu$ as a function on the symmetric group,
or equivalently as an element of its group algebra.
Then, denoting by $\sigma_j$ the elementary transposition
$(j,j+1)$,
the induction $\nu \rightarrow \mu= \nu \sigma_j$, 
implies the recursion
\begin{equation}\label{RecY}
A_\mu = A_\nu Y_j(u_{\nu(j)} ,u_{\nu(j+1)} ) \ , 
\end{equation}
where 
$\displaystyle Y_j(u,v) = { {\rm i} (u-v)\sigma_j +c  \over {\rm i}(u-v) -c  }$
satisfies the Yang-Baxter equation 
\begin{equation}\label{YB}
Y_i(u,v)Y_{i+1}(u,w)Y_i(v,w) =
Y_{i+1}(v,w)Y_i(u,w)Y_{i+1}(u,v) \ .
\end{equation}
The parameters $u_i$ are obtained as solutions of a system of
nonlinear algebraic equations. The next problem
is to expand the operator $A_\mu$ on the basis of permutations. Yang
obtained this expansion by means of
the recurrence relation (\ref{RecY}), and as far as we know, no closed
expression for the coefficients $A_\mu(\tau)$ is known.
One aim of this note is to investigate this question. We give a partial
answer to the original problem, and a complete one to some variants
that we explain below.

Getting rid of the normalization constants, we replace Yang's operators by 
$$
Y_i(u,v)=1-(u-v)\sigma_i \ .
$$
More general solutions of (\ref{YB})
can be obtained by replacing
the elementary transpositions  $\sigma_i$ by the generators $T_i$ of
a Hecke algebra, satisfying the braid relations and a quadratic relation
$$
(T_i-q_1)(T_i-q_2)=0\ . 
$$
For generic values of $q_1,q_2$ this solution is
$$
Y_j=   1 + {v/u -1 \over q_1 +q_2  } T_j  \ , 
$$
while for instance, in the degenerate case $T_i^2=0$ one can take $Y_i=1-(u-v)T_i$.

To each of these solutions and to each system of spectral parameters, one
can associate a basis $\{Y_\mu(\u)\, |\, \mu\in\S_n\}$
of the corresponding Hecke algebra, called the Yang-Baxter basis.
From this data, we define a bilinear form on the Hecke algebra. 
A fundamental property of the Yang-Baxter basis is 
that it is essentially self-adjoint with respect to
this form ({\it i.e.}
its adjoint basis is obtained by  permuting  the indices and rescaling).
We also give closed formulas for the coefficients of the Yang-Baxter
basis on the usual basis $\{T_\mu\}$ in the two degenerate cases
$T_i^2=0$ and $T_i^2=-T_i$. 
The solution is obtained in terms of specializations of double
Schubert and Grothendieck polynomials, which were originally
defined as canonical bases of the cohomology and Grothendieck
rings of flag manifolds.
The first case can be interpreted as
giving the terms of lowest degree in Yang's coefficients.
The existence of a connection between Schubert polynomials and
the Yang-Baxter equation was first noticed by Fomin and Kirillov \cite{FK1}.
All that we present here is for the root system $A_n$. Generalizations
of Schubert polynomials for types $B_n,C_n,D_n$ are given in \cite{PR,Ful}
(see also \cite{P}).
The case $B_n$ is studied in detail in \cite{FK2} and \cite{BH}.

This paper is organized as follows. In sections \ref{FLAGS} and \ref{SCHUBGRO}
we recall the definitions of Schubert and Grothendieck polynomials, and
their geometric interpretations. Then we introduce the Yang-Baxter basis
(Section \ref{YBbasis}) and the bilinear form on the Hecke algebra
with respect to which this basis is well behaved (Section \ref{ORTHO}).
Finally, we present the expansions of the YB basis on the standard basis
in terms of Schubert and Grothendieck polynomials (Section \ref{EXP}).

\section{Flag varieties and Hecke algebras}\label{FLAGS}

Hecke algebras arise naturally in the theory of flag manifolds,
as algebras of operators acting on the cohomology \cite{BGG,Dem}
or on the Grothendieck ring \cite{Dem,Lus}.
Recall that a flag of vector spaces is an  increasing 
sequence of vector spaces 
$V_1 \subset V_2 \subset\cdots \subset {\C}^n$. A flag is said to be
 complete when it is of the type
$$V_1 \subset V_2 \subset\cdots \subset V_n= {\C}^n  \ ,$$
with ${\rm dim\,}(V_i) = i$, $i=1,\ldots, n$.
The set of complete flags in ${\C}^n$ is an algebraic projective variety
${\cal F}_n={\cal{F}}({\C}^n)$. 
This variety is equipped with  tautological line bundles
$L_1$, $L_2, \ldots, L_n$ which are defined as follows:
for each $1\leq i\leq n$, the collection 
$W_i := \{(F,V_i) ,\, V_i\in F, F\in {\cal{F}}({\C}^n) \}$ 
is a vector bundle on ${\cal{F}}({\C}^n)$. 
Then, writing $E^\vee$ for the dual
of a vector bundle $E$, 
$$L_i := (W_i /\, W_{i-1} )^\vee  \ .$$
The  cohomology and Grothendieck ring (of classes of vector bundles) 
of the flag variety are both quotients of the ring of polynomials in
$n$ variables. More precisely, writing $1+\gamma_i$ for the Chern class
of $L_i$, and $\xi_i$ for the class of $L_i$ in the Grothendieck ring,
 $1\leq i\leq n$, the cohomology ring is 
$$\C[\gamma_1,\ldots,\gamma_n]=
\C[x_1, \ldots, x_n] / {\cal{J}}(0)$$
and the Grothendieck ring
$$\C[\xi_1,\ldots,\xi_n] = 
\C[x_1, \ldots, x_n] / {\cal{J}}(1) \ ,$$
where ${\cal{J}}(0)$ is the ideal generated by the 
$f(x_1,\ldots, x_n) - f(0,\ldots, 0)$, $f$ symmetrical,
and ${\cal{J}}(1)$ is the ideal generated by the 
$f(x_1,\ldots, x_n) - f(1,\ldots, 1)$.

Both rings are graded modules of rank $n!$
with natural bases
parametrized by the symmetric group  $\S_n$. The Poincar\'e polynomial
of both spaces is 
$$ (1+q)\,  (1+q+q^2)\cdots  (1+q+\cdots + q^{n-1} ) \ .$$

The equivariant Grothendieck ring $K_{GL(n)}({\cal F}_n)$ is the ring
of Laurent polynomials $\Z [x_1^{\pm 1},\ldots, x_n^{\pm 1}]$,
where $x_i$ is now the equivariant class of $L_i$.
For the purpose of introducing
bases depending on parameters, we shall work in the extended ring
$R:= \C[y_1,\ldots,y_n;x_1^{\pm 1},\ldots, x_n^{\pm 1}]$.
This ring is a free module of rank $n!$ over $\Sym(X)[t]$, where
$t=(x_1\cdots x_n)^{-1}$ and $\Sym(X)$ is the ring of symmetric
polynomials in the $x_i$'s with coefficients in $\C[y_1,\ldots,y_n]$.

The flag variety can be constructed as an iterated sequence
of projective bundles (see \cite{H}). This construction gives rise to
various operators on the cohomology and Grothendieck rings, corresponding
to natural operations on the ring of polynomials:

implies the existence
of natural linear operators acting on the ring of polynomials:

\begin{itemize}
\item permutations, generated by the elementary transpositions
$$
\sigma_i \ : \ \ f\longmapsto \sigma_i f = f(\ldots,x_{i+1},x_i,\ldots)
$$
\item divided differences
$$
\partial_i f = {f-\sigma_if\over x_i-x_{i+1}}
$$
used in \cite{BGG} and \cite{Dem} to construct a basis of the cohomology,
\item $s_i =\sigma_i+\partial_i$ (see \cite{LLT}),
\item isobaric divided differences $\pi_i(f)=\partial_i(x_if)$,
which appear in Demazure's character formula \cite{DemCar},
\item $\bar\pi_i =\pi_i-1$  ({\it cf.} \cite{Lus}).
\item $T_i =-(q_1+q_2) \bar\pi_i +q_2 \sigma_i$ ({\it cf.} \cite{DKLLST}).
\end{itemize}

Each of these families of operators  satisfies the braid relations
\begin{equation}
D_i D_{i+1} D_i = D_i D_{i+1} D_i \ , 
\end{equation}
\begin{equation}
D_i D_j = D_j D_i \quad |i-j|>1\ . 
\end{equation}
This ensures the existence of operators $D_\mu$, $\mu\in  \S_n$, which are 
obtained by taking, for a given permutation $\mu$, an arbitrary 
reduced decomposition  $\mu= \sigma_{i_1} \sigma_{i_2}\cdots\sigma_{i_r}$,
and putting $D_\mu = D_{i_1} D_{i_2} \cdots D_{i_r}$

The simple operators $D_i$  are of the Hecke type, {\it i.e.}
satisfy a quadratic relation:
\begin{equation}
 \sigma_i^2=1, \quad \partial_i^2 =0, \quad s_i^2=1, \quad \pi_i^2 = \pi_i,
\quad \overline{\pi}_i^2 = - \overline{\pi}_i ,
\quad (T_i-q_1)(T_i-q_2) = 0 \ . 
\end{equation}

\vskip 3mm
We shall write  ${\cal{H}}^\sigma$, ${\cal{H}}^\partial$, ${\cal{H}}^s$, 
${\cal{H}}^\pi$, ${\cal{H}}^{\bar \pi}$ and ${\cal{H}}^T$ for 
the different algebras
generated by the corresponding simple operators. 
The algebra generated by the $T_i$'s together with the variables $x_i$ 
(considered as operators $f\mapsto x_if$) is the affine Hecke algebra
\cite{Ch1,Ch2}. The other algebras are specializations of this
one. In particular, the semidirect product ${\cal H}^s\times \C[x_1,\ldots,x_n]$
is the degenerate affine Hecke algebra \cite{Ch,Ch1}. 
Note that  we distinguish between
${\cal{H}}^\sigma$ and ${\cal{H}}^s$ which have a different
action on $R$.

\section{Schubert polynomials and Grothendieck polynomials}\label{SCHUBGRO}

The ring  $R$ admits two  distinguished bases as a free module
over the ring $R^{\S(X)}$ of symmetric polynomials in $X$, 
the {\it double Schubert polynomials} $X_\mu$
and the {\it double Grothendieck polynomials} $G_\mu$, $\mu\in \S_n$,
which are defined as follows.
For $\mu = \omega := (n,\ldots, 2,1)$, 
$$X_\omega := \prod_{i+j\leq n} (x_i -y_j)\quad , \quad 
G_\omega := \prod_{i+j\leq n} (1- y_j/x_i )   $$
and otherwise
\begin{equation}\label{DOUBLE}
X_\mu:= \partial_{\mu^{-1}\omega}X_\omega \,
\end{equation}
\begin{equation}
G_\mu:=\pi_{\mu^{-1}\omega}G_\omega \ .
\end{equation}

The specialization $X_\mu( y_1=0,\ldots,\,  y_n=0)$ is a representative
of a Schubert cycle in the cohomology \cite{LS,Mcd}.
The specialization
$G_\mu( y_1=1, \ldots ,\,  y_n=1)$ is a representative of
the class of the structure sheaf of a Schubert variety  in the
Grothendieck ring \cite{FL}.

Double Schubert polynomials also have an interpretation as universal
polynomials for degeneracy loci ({\it cf.} \cite{FulXX}).
The simplest case corresponds to a pair of vector bundles $E,F$, and a map
$f: E\longrightarrow F$. The $r$-th degeneracy locus of this map is the set
of points where the corank of $f$ is $\ge r$. 
Under suitable genericity conditions, the class of this
degeneracy locus is a polynomial in the Chern classes of $E$ and
$F$, which can be identified with some Schubert polynomial in the
Chern roots of $E$ and $F$. General Schubert polynomials correspond
to a pair of flags of vector bundles, a map $f$ between them,
and corank conditions on $f$.

\section{Yang-Baxter bases of Hecke algebras}\label{YBbasis}

We first define the elementary Yang-Baxter operators $Y_i$
for the various realizations of Hecke algebras considered in
Section \ref{FLAGS}.
The Yang-Baxter basis will then be constituted by the operators
$Y_\mu$ ($\mu\in\S_n$), defined inductively by
$$
Y_\mu(\u)= Y_\nu(\u) Y_j(u_{\nu(j)} ,u_{\nu(j+1)} )
$$
for $\mu=\nu\sigma_j$, $\ell(\mu)>\ell(\nu)$.
The elementary operators are as follows:
$$Y_j^\sigma(u,v) := 1 - (u-v) \sigma_j  \ \in {\cal{H}}^\sigma  $$
$$ Y_j^\partial(u,v) :=  1 - (u-v) \partial_j \ \in {\cal{H}}^\partial \ ,  $$
$$Y_j^{\bar \pi}(u,v):=  1+ (1-v/u)  \bar \pi_j \ \in {\cal{H}}^{\bar \pi}_ , $$
$$Y_j^T(u,v;q_1,q_2) :=   1 + {v/u -1 \over q_1 +q_2  } T_j \ \in {\cal{H}}^T \ . $$
One can recover $Y_j^\partial$, $Y_j^{\bar\pi}$ and $Y_j^\sigma$
from the last case $Y_j^T$ by various  specializations. 
Thus, setting 
$$
\tilde{Y}_j^T(u,v;q_1,q_2) =
Y_j^T\left( e^{(q_1+q_2)u},\, e^{(q_1,+q_2)v)}; q_1,q_2 \right) \ ,
$$
one has
$$
Y_j^\sigma(u,v)= \lim_{q_1+q_2\rightarrow 0} \tilde{Y}_j^T(u,v;q_1,q_2) \ ,
$$
$$
Y_j^{\bar\pi}(u,v)=\tilde{Y}_j^T(u,v;-1,0) \ .
$$
The operator $Y_j^\partial$ can also be obtained from
$Y_j^T$ by combining a similar specialization with an appropriate
homographic substitution on the variables $x_i$.

In summary, given an arbitrary choice of the parameters $\u = (u_1,\ldots, u_n)$, 
for all the above Hecke algebras  
there exists a linear basis $\{Y_\mu(\u), \, \mu\in\S_n \}$, that we
call the {\it Yang-Baxter basis}.
The problem
that we shall examine is to express the Yang-Baxter basis in the standard 
bases $\{\mu \}$, $\{\partial_\mu \}$, $\{\bar \pi_\mu \}$ or $\{T_\mu \}$.

We shall deduce all relations from the following recursion 
(in the case of the generic Hecke algebra)
\begin{equation}
Y^{T}_{\mu \sigma_j}(\u) = Y^{T}_\mu(\u) \cdot (1+ {u-1\over q_1+q_2} T_j )\ , 
\end{equation}
when $\ell(\mu \sigma_j)>\ell(\mu)$, 
with $u = u_{\mu(j+1)}/ u_{\mu(j)}$.
We note that, in the case $\ell(\mu \sigma_j)<\ell(\mu)$, one has, 
with the same $u$, 
\begin{equation}
Y^{T}_{\mu \sigma_j}(\u)  \cdot (1+ {1/u -1\over q_1+q_2} T_j )
 = Y^{T}_\mu(\u) \cdot (1- {2-u-1/u\over (q_1+q_2)^2} q_1q_2 ) \ .
\end{equation}
More generally, the Hecke relation implies, for all $j,u,v $,
\begin{equation}\label{INDU}
(1+ {u -1\over q_1+q_2} T_j ) \cdot 
(1+ {v -1\over q_1+q_2} T_j )
=  1- {(u -1)(v-1)q_1q_2\over (q_1+q_2)^2}
+ {uv-1\over q_1+q_2 }   T_j   \ .
\end{equation}
This product is thus a constant iff $u = 1/v$.

Cherednik \cite{CheDuke} has shown that one can recover orthogonal idempotents
in the group algebra of the symmetric group by some specialization of the
Yang-Baxter basis. 
He uses them to describe bases of
representations of the symmetric group corresponding to skew partitions.

\section{Orthogonality properties of the Yang-Baxter basis}\label{ORTHO}

We shall only treat the case of the generic Hecke algebra, and 
recover the other ones by specialization.

Define an anti-automorphism $\varphi $ on ${\cal{H}}[\u]$ by 
\begin{equation}
\varphi(T_\mu)=T_{\mu^{-1}},\qquad \varphi(u_i)=u_{n-i+1}
\end{equation}
and a bilinear form $<\, ,\, >$  by
\begin{equation}
<h_1\, ,\, h_2> := h_1\cdot \varphi(h_2)\big|_{T_\omega}  \ , 
\end{equation}
for $h_1,h_2 \in{\cal{H}} $
where $ h\, |_{T_\omega}$ denotes the coefficient of $T_\omega$
in the expansion of $h\in{\cal{H}}$ in the basis $T_\mu$. 

\begin{theorem} \label{ORTHOG}
Let $\u$ be an arbitrary system of spectral parameters, and
$$
\Delta^T(\u) := \prod_{1\leq i<j\leq n} {u_j/u_i - 1\over q_1+q_2}
$$ 
Then the Yang-Baxter basis $\{Y_\mu^T \}$, $\mu\in\S_n$,
is adjoint to $\{  Y_{\omega\mu}^T\,  \Delta^T(\u^{\omega\mu})^{-1} \}$, 
 {\em i.e.} one has
$$<Y_\mu^T\, ,\, Y_\nu^T> =  \Delta^T(\u^{\omega\mu}) \delta_{\nu,\omega\mu}
$$
\end{theorem}

\Proof    Let $\mu$ and $j$ be such that $\ell(\mu\sigma_j)>\ell(\mu)$. Then
$$
Y_\mu^T\cdot (1 + {u-1 \over q_1+q_2} T_j) \cdot \varphi(Y_\nu^T)
=
Y_\mu^T\cdot \varphi \bigl( Y_\nu^T \cdot ( 1 + {1/u -1 \over q_1+q_2} T_j) \bigr)  
$$
Now, by induction on the length of $\mu$, we can suppose that we know all
the pairings $<Y_\mu^T\, , \, Y_\eta^T>$. Since for any constant $c$, 
$Y_\nu^T (1+ c T_j)$ is a linear combination of
$T_\nu$ and $T_{\nu \sigma_j}$, the  non-zero pairings
can occur only  for $\mu= \nu $ or $\mu=\nu\sigma_j$. 
The conclusion follows from (\ref{INDU}). \cqfd

\begin{corollary} 
Given  an arbitrary set of parameters   $\u$,
any element $h $ of the Hecke algebra ${\cal{H}}^T $ can be expressed as
$$  h= \sum_{\mu}  {1 \over \Delta^T(\u^{\omega\mu})   } 
<h\, ,\, Y_{\omega\mu}^T(\u) >  Y_\mu^T({\bf u}) \ . $$
\end{corollary}

\vskip 3mm
\begin{example}
{\rm
We take   $n=3$, and we write for short $(ji)T_k$ instead of 
$ (1+ {u_j/u_i -1 \over q_1+q_2} T_k)   $ and $Y_\mu$ instead of
$Y_\mu^T$. Then, the
Yang-Baxter basis is

$$\matrix{
     &Y_{123}=1& \cr
Y_{213}= (21)T_1  &   & (32)T_2=Y_{132}  \cr
Y_{231}= (21)T_1 (31)T_2  &   & (32)T_2 (31)T_1 = Y_{312}  \cr
(21)T_1 (31)T_2 (32)T_1  &=Y_{321}=   & (32)T_2 (31)T_1 (21) T_2  \cr }$$
and its image under $\varphi$ is
$$\matrix{
     &1 \cr
(23)T_1  &   & (12)T_2  \cr
(13)T_2 (23)T_1  &   & (13)T_1 (12)T_2  \cr
(23)T_1 (13)T_2 (12)T_1  &=   & (12)T_2 (13)T_1 (23) T_2  \cr }$$

Conversely, the expressions of the elements  $T_\mu$ 
in the Yang-Baxter basis are, writing $\Delta$ for $\Delta^T(\u)$:

$\Delta\,  T_{123}= \Delta Y_{123}$; 
$\Delta\, T_{213}= {(u_3/u_2-1) (u_3/u_1-1) \over (q_1+q_2)^2  }
( Y_{213} -Y_{123})   $

$\Delta\, T_{132}= {(u_2/u_1- 1) (u_3/u_1 -1) \over (q_1+q_2)^2  }
( Y_{132} -Y_{123})   $

$\Delta\, T_{231}=  { u_3/u_2 -1 \over q_1+q_2} (Y_{231}- Y_{213})
  - { u_3/u_1 -1 \over q_1+q_2} (Y_{132}- Y_{123})   $

$\Delta\, T_{312}=  { u_2/u_1 -1 \over q_1+q_2} (Y_{312}- Y_{132})
  - { u_3/u_1 -1 \over q_1+q_2} (Y_{213}- Y_{123})   $

$ \Delta\,  T_{321}
= Y_{321} - Y_{231} - 
Y_{312} + Y_{213}  +Y_{132}
 + ( 1 - {1+u_3/u_1 -u_3/u_2 -u_2/u_1\over 
(1+q_1/q_2)(1+q_2/q_1)   }      ) Y_{123}  $
}
\end{example}

\vskip 3mm
\section{Explicit coefficients}\label{EXP}

An easy way to get the sequence of parameters used in the 
factorized expression of a Yang-Baxter element is
provided by a planar representation of a permutation due to
Rothe (1800), that we shall call its Rothe-Riguet diagram \cite{Mcd}.

Represent a permutation $\mu $, in the Cartesian plane $\N\times\N$,
by the set of points of coordinates $(i,\mu(i))$. The
diagram of $\mu$ is the set of points (denoted by boxes $\carre$)
$\{ (i, \mu(j) ) \, :\, i<j, \mu(i)>\mu(j)    \}$. This diagram is a convenient
way to represent the inversions of a permutation. The number of boxes
in the diagram of $\mu$ is exactly the length $\ell(\mu) $ of $\mu$.
Each box  $(i, \mu(j))$
is then filled with the factor
$$1+ { u_{\mu(j)}/u_{\mu(i)} -1\over q_1+q_2  } T_k  \, $$
where $k- i $ is the number of boxes 
in the same column strictly below the given box.
Thus in each column, the indices of the $T_k$ are increasing by 1
upwards, starting from the index of the column.

\begin{example}\label{EX}{\rm
Let $\mu= 31542$. Still writing $(ji)T_k$ for 
$(u_j/u_i -1)/(q_1+q_2)\, T_k$, and denoting the points
$(i,\, \mu(i))$ by a $\star$, the diagram of $\mu$ is

\def\s{\scriptstyle}
$$\matrix{
\s{5} & &\cdot    &\star  & \cdot&\cdot&\cdot\cr
\s{4}&&\cdot      &(54)T_4& \cdot&\star&\cdot \cr
\s{3}&&\star      &\cdot  & \cdot&\cdot&\cdot   \cr
\s{2}&&(32)T_2    &(52)T_3& \cdot&(42)T_4&\star\cr
\s{1}&&(31)T_1    &(51)T_2& \star&\cdot   &\cdot\cr \cr
\mu   &=  &\s{3}&\s{5}&\s{1}&\s{4}&\s{2}       \cr
}$$

The following proposition, obtained by induction on the length $\ell(\mu)$
tells us that 
$$Y_\mu^T= (54)T_4 \, (32)T_2 \, (52)T_3 \, (42)T_4 \, (31)T_1 \, (51)T_2  $$
}
\end{example}

\vskip 3mm
\begin{proposition}
The reading of the Rothe diagram of a permutation $\mu$
(in the occidental way, from left to right and top to bottom), 
as the product of the factors contained in the boxes, 
gives an expression of the Yang-Baxter element $Y_\mu^T(\u)$. \cqfd
\end{proposition}

Let us now consider degenerate Hecke algebras.
The following formulas could in principle be derived from
Theorem \ref{ORTHOG}, but it will be more convenient to provide direct
proofs in each case.
We first treat the case of the simplest Hecke relation
$D^2_i =0$, {\it i.e.} the algebra ${\cal{H}}^\partial$.

\smallskip 
\begin{proposition}
The entries of the transition matrix   
$Y_\mu^\partial(\u) \rightarrow \partial_\nu$ of ${\cal H}^\partial$
are specializations of double Schubert polynomials. 
More precisely,
\begin{equation}\label{SCHUB}
Y_\mu^\partial(\u) = \sum_\nu X_\nu (\u^\mu,\u)\,  \partial_\nu \ .
\end{equation}
\end{proposition}

\Proof    We use  the induction $\mu \rightarrow \mu\sigma_k$ 
with $\ell(\mu\sigma_k)>\ell(\mu)$.
Assume that (\ref{SCHUB}) holds for $\mu$, \ie that
$$
Y_\mu^\partial = \sum_{\nu :\ell(\nu\sigma_k) = \ell(\nu)+1} 
\left[ 
X_\nu(\u^\mu,\u) \partial_\nu 
+ X_{\nu\sigma_k}(\u^\mu,\u) \partial_{\nu\sigma_k}
\right]
$$
Multipliying by the factor 
$(1 - (u_i-u_j) \partial_k )$, where $i= \mu(k)$, $j=\mu(k+1)$,
one obtains the coefficients of $\partial_\nu$ and $\partial_{\nu\sigma_k}$ in 
$Y_{\mu\sigma_k}$. The one of $\partial_{\nu\sigma_k}$ is equal to
$X_{\nu\sigma_k}({\bf u}^\mu,{\bf u})-(u_i-u_j)X_\nu({\bf u}^\mu,{\bf u})$. 
Taking into account the relation
$$(u_i -u_j) X_{\nu}(\u^\mu, \u) =
X_{\nu\sigma_k}(\u^\mu, \u) - X_{\nu\sigma_k}(\u^{\mu\sigma_k} , \u)    $$
which follows directly from the definition (\ref{DOUBLE}),
we see that this coefficient is indeed given by (\ref{SCHUB}).  \cqfd


There are several related formulas involving Schubert polynomials and
the nilCoxeter algebra. The first one is a recursion on Schubert
polynomials  (\cite{LSschub2}, Prop. 1.7 and 1.12,
see also \cite{Mcd}, (7.2) p. 95) in ${\cal H}^\partial$ (called
{\it alg\`ebre des diff\'erences divis\'ees} in \cite{LSschub2,LS}).
This recursion is directly equivalent to the generating function
obtained in \cite{FS,FK1} by a different method involving the Yang-Baxter
equation. The coefficients $X_\nu({\bf u}^\mu,{\bf u})$ also
appear in the expansion of permutations in the divided differences
algebra ${\cal H}^\partial$. According to \cite{LS,LSdd}, one has
\begin{equation}\label{decperm}
\mu = \sum_\nu X_\nu({\bf x}^\mu,{\bf x})\partial_\nu  \ .
\end{equation}
This follows from the Newton interpolation formula, since for
any polynomial $f$,
$$
\sum_\nu X_\nu({\bf x}^\mu,{\bf y})\partial_\nu f({\bf y}) = f({\bf x}^\mu) \ .
$$
Note however that in (\ref{decperm}) the divided differences
$\partial_\nu$ act on the variables $x_i$. To recover permutations
from the Yang-Baxter operators, one just has to compose them
with the specialization homomorphism $\xi : u_i\mapsto x_i$. In other
terms, if one defines normal ordering  in the divided differences
algebra by the condition that for a monomial $M$ in the $x_i$
and the $\partial_j$, $:\! M\!:$ is the same expression where all the
$x_i$ have been moved to the left of the difference operators,
without alterating the relative order of the $\partial_j$, one has

\begin{corollary}
As an operator on polynomials,
the permutation $\mu$ is equal to the normal ordering of the Yang-Baxter
element $Y_\mu^\partial({\bf x})$
$$
\mu =\ :\!Y_\mu^\partial ({\bf x})\!:   \ =\ \xi\circ Y_\mu^\partial({\bf u}) \ .
$$
\end{corollary}
For example,
$$
Y_{321}^\partial ({\bf u}) =
(1-(u_2-u_3)\partial_1)(1-(u_1-u_3)\partial_2)(1-(u_1-u_2)\partial_1)
$$
$$
= 1+\cdots + (u_3-u_2)(u_3-u_1)(u_2-u_1)\partial_{321}
$$
and, writing $\sigma_i=1-(x_i-x_{i+1})\partial_i$,
$$
(321)=\sigma_1\sigma_2\sigma_1
=
(1-(x_1-x_2)\partial_1)(1-(x_2-x_3)\partial_2)(1-(x_1-x_2)\partial_1)
$$
$$
= 1+\cdots +(x_2-x_1)(x_3-x_1)(x_3-x_2)\partial_{321} = \ :Y^\partial_{321}({\bf x}):
$$
the reduction to a normally ordered expression involving
repetitive use of the Leibniz formula
$\partial_i(PQ)=(\partial_iP)Q+\sigma_i(P)\partial_iQ$.

One could think of obtaining the expansion (\ref{SCHUB}) by taking
the specialization of the generating function of double Schubert polynomials.
For example, for $\mu=(1324)$, the generating function of Fomin and
Kirillov specializes to
$$
(1+(x_2-x_1)\partial_3)(1+(x_3-x_1)\partial_2)(1+(x_1-x_1)\partial_1)
(1+(x_3-x_2)\partial_3)(1+(x_1-x_2)\partial_2)(1+(x_1-x_3)\partial_3)
$$
and it is not immediate that this expression is equal to
the Yang-Baxter element
$$
Y^\partial_{1324} = 1-(x_2-x_3)\partial_2 \ .
$$
In other words,  formula (\ref{SCHUB}) implies
\begin{corollary}
Let ${\cal F}({\bf x},{\bf y})$ be the Fomin-Kirillov generating
function for double Schubert polynomials (see \cite{FK1}). Then,
$$
{\cal F}({\bf u}^\mu,{\bf u}) = Y^\partial_\mu({\bf u}) \ .
$$
\end{corollary}

\vskip 3mm
Next, we  turn back to the original case of Yang and Baxter,
{\it i.e.} the case where
${\cal{H}} = {\cal{H}}^\sigma$ is the group algebra of the symmetric group.

\begin{proposition}
The  Yang-Baxter coefficients $A_\nu(\mu)$, 
which are the coefficients of the expression of the 
Yang-Baxter basis in the basis of permutations,
are non-homogeneous polynomials whose leading terms
(i.e. term of lowest degree) are equal to the specializations  
$X_\nu(\u^\mu,\u)$ of double Schubert polynomials.
\end{proposition}

\Proof  In the case of ${\cal{H}}^\partial$, expanding a Yang-Baxter element 
$Y_\mu^\partial(\u)$,
after having chosen a reduced decomposition  $w$ of $\mu$, 
involves only subwords of $w$  which are reduced decompositions
(because products of $\partial_i$ which are not reduced decompositions
vanish).
These subwords give exactly the 
same contribution in the case of ${\cal{H}}^\sigma$. Extra terms 
involve at least  a square $ \cdots (u_i-u_j) \sigma_k (u_{i'}- u_{j'}) \sigma_k 
\cdots$ and thus give a contribution to $A_\nu(\mu)$ which is of
degree strictly bigger than $\ell(\nu)$. 

For example, for $\S_3$, 
the only coefficient which does not coincide with a Schubert polynomial
is $A_{123}(321)   = 1 + (u_1-u_2) (u_2-u_3)$. \cqfd

\vskip 3mm
Another interesting case for geometry is the algebra ${\cal{H}}^{\bar \pi}$.

\begin{proposition}
The entries of the transition matrix  
 $Y_\mu^{\bar \pi}(\u) \rightarrow \bar \pi_\nu$ of ${\cal{H}}^{\bar \pi}$
are specializations of double Grothendieck polynomials. 
More precisely,
$$Y_\mu^{\bar \pi}(\u) = \sum_\nu G_\nu (\u^\mu,\u)\,  \partial_\nu \ .$$
\end{proposition}

The proof is the same as in the case of the $Y^\partial_\mu(\u)$,
the relations $\partial_i^2=0$ being replaced by $\pi_i \bar \pi_i =0$.
Remark that the expansion of a product of $k$ factors
$1+ (u-1)\bar \pi_i$ involves $2^k$ terms, in contrast with the case
of $1+ (u-1)\partial_i$ where only those terms which are 
reduced decompositions give a non zero contribution.


Fomin and Kirillov \cite{FK} have given a generating function
for Grothendieck polynomials in the algebra ${\cal H}^{\bar\pi}$,
related to the recursion of \cite{LSschub2} (Th. 2.5) and to formula
(4.11) of \cite{LSdd}. One has the same type of specialization
as for Schubert polynomials:
\begin{corollary}
Let ${\cal G}$ be the Fomin-Kirillov generating function of double
Grothendieck polynomials. Then,
$$
Y_\mu^{\bar\pi}(u_1^{-1},\ldots,u_n^{-1})
= {\cal G}({\bf u}^\mu,{\bf u}) \ .
$$
\end{corollary}


One can recover Schubert polynomials from Grothendieck polynomials.
Indeed, the change of variables $x_i \rightarrow 1/(1-a_i)$,
$y_j \rightarrow 1/(1-b_j)$ transforms the maximal
Grothendieck polynomial $G_\omega(X,Y)$ into 
$\prod_{i,j: i+j\leq n} (a_i -b_j)/(1-b_j)  $, i.e. into the Schubert
polynomial $X_\omega(A,B) $    , up to a factor
independent of the $a_i$.  On another hand, 
the operator $\pi_i$ becomes 
$f \longrightarrow \partial_i (  (1 + a_{i+1})\cdot f) $, where 
the $\partial_i$ are now the operators corresponding to the variables 
$\{ a_i \}$. Thus, the term of smallest degree (in the $a_i$)
 of the polynomial
$G_\mu(x_i,y_j)$, when expressed in terms of the $a_i$, $b_j$,
is the Schubert polynomial $X_\mu(a_i, b_j)$.

\section{Appendix}

\subsection{Action on the ring of polynomials}

We have already mentioned that the different Hecke algebras naturally
occur as algebras of operators on the ring $\C[x_1,\ldots,x_n]$.
This action is in fact related to the $\chi_y$-characteristic
of Hirzebruch \cite{H} (see \cite{DKLLST}).
However, in the preceding sections we did not make use of the action
on the variables $x_i$. Let us  mention only two applications \cite{DKLLST,LLT}.

In the case where the spectral parameters are chosen to be
$u_i=q^{i-1}$, with $q_1=q,\, q_2=-1$, the Yang-Baxter element
$Y_\mu$, for $\mu$ the maximal element of a Young subgroup
$\S(I)=\S_{i_1}\times\S_{i_2}\times\cdots\times\S_{i_r}$
factorizes into
$$
Y_\mu (f)= \Delta_1 \Delta_2 \cdots \Delta_r \, \partial_\mu (f)
$$
where $\Delta_k$ is the $q$-Vandermonde associated to the $k$-th
factor of the Young subgroup, {\it i.e.}, setting $m_k=i_1+\cdots +i_k$
$$
\Delta_k 
=
\prod_{m_k+1 \le i < j \le m_{k+1}} (qx_j-x_i) \ .
$$
In \cite{DKLLST} is described how to use these special
elements to obtain bases of the irreducible representations
of the Hecke algebra, as well as $q$-idempotents generalizing
the Young idempotents.

The case of ${\cal H}^s$, where the parameters are now
$u_i=i$, is given in \cite{LLT}.
The images of the monomial $x^{n-1}_1 x_2^{n-2}\cdots x_1^0$
under Yang-Baxter elements $Y_\mu^s$ constitute a linear basis of
the cohomology of the flag manifold, which can be regarded as
a deformation of the Schubert basis.

Yang-Baxter elements corresponding to the maximal elements of a
Young subgroup still satisfy a factorization property. 
For example, if $\omega$ is the permutation $(n,n-1,\ldots,1)$,
one has
$$
Y_\omega^s (f) = \prod_{1\le i<j \le n}
(1+x_j-x_i) \, \partial_\omega (f) \ .
$$
The polynomial $ \prod_{1\le i<j \le n} (1+x_j-x_i) $ can be
interpreted as the total Chern class of the flag variety.

\subsection{Examples}

For $n=4$, the double Schubert polynomials are:
\small

\medskip \noindent
$ X_{4321} = (x_1-y_1)(x_1-y_2)(x_1-y_3) (x_2-y_1)(x_2-y_2) (x_3-y_1)$

\medskip\noindent
$X_{3421} = X_{4321}/(x_1-y_3)$ ; 
$X_{4231} = X_{4321}/(x_2-y_2)$ ; $X_{4312} = X_{4321}/(x_3-y_1)$

\medskip \noindent
$X_{2431} = (x_1-y_1)(x_2-y_1)(x_3-y_1)(x_1+x_2-y_2-y_3)$; 
$X_{3241} = X_{4321}/(x_1-y_3)(x_2-y_2)$;
$X_{3412} = X_{4321}/(x_1-y_3)(x_3-y_1)$;
 $X_{4132} = (x_1-y_1)(x_1-y_2)(x_1-y_3)(x_2+x_3-y_1-y_2)$;
$X_{4213} = X_{4321}/(x_2-y_2)(x_3-y_1)$

\medskip \noindent
$X_{2341} = (x_1-y_1)(x_2-y_1)(x_3-y_1)$;
 $X_{1432} = x_1^2x_2+x_1^2x_3+x_1x_2^2+x_1x_2x_3 +x_2^2x_3 
-(x_1^2+x_1x_2+x_2^2)(y_1+y_2) - (x_1x_2+x_1x_3+x_2x_3)(y_1+y_2+y_3) 
+(x_1+x_2)(y_1^2+y_1y_2+y_2^2) 
+(x_1+x_2+x_3)(y_1y_2+y_1y_3+y_2y_3) 
-(y_1^2y_2+y_1^2y_3+y_1y_2^2+y_1y_2y_3+y_2^2y_3) $; 
$X_{2413} = 
(x_1-y_1)(x_2-y_1)(x_1+x_2-y_2-y_3)$; 
$X_{3142} = (x_1-y_1)(x_1-y_2)(x_2+x_3-y_1-y_2)$;
$X_{3214} = (x_1-y_1)(x_1-y_2)(x_2-y_1)$; 
$X_{4123} = (x_1-y_1)(x_1-y_2)(x_1-y_3)$.

\medskip \noindent
$X_{1342} = x_1x_2+x_1x_3+x_2x_3 -(y_1+y_2)(x_1+x_2+x_3)+y_1^2 +y_1y_2+y_2^2$ ;
$X_{1423} =x_1^2 +x_1x_2+x_2^2 -(x_1+x_2)(y_1+y_2+y_3)+y_1y_2+y_1y_3+y_2y_3 $ ;
$X_{2143} =(x_1-y_1)(x_1+x_2+x_3-y_1-y_2-y_3)$;
$X_{2314} = (x_1-y_1)(x_2-y_1)$ ;
$X_{3124} = (x_1-y_1)(x_1-y_2)$ .

\medskip \noindent
$X_{1243} = x_1+x_2+x_3-y_1-y_2-y_3$ ;
$X_{1324} = x_1+x_2-y_1-y_2$ ; $X_{2134} = x_1-y_1$ .

\medskip\noindent
$X_{1234} = 1$

\medskip\normalsize
\goodbreak
The Grothendieck polynomials for $n=3$ are

\setbox21=\hbox{$G_{321}=(1-{y_1 \over x_1 })$}
\setbox22=\hbox{\kern 1cm $(1-{y_1 \over x_2 })(1-{y_2 \over x_1 }) $}
\setbox2=\vbox{\box21 \box22}
\setbox31=\hbox{$ \pi_1 $}
\setbox32=\hbox{$\, \swarrow$}
\setbox3=\hbox{\kern 25mm{\raise 3mm\box31}\box32}
\setbox41=\hbox{$ \searrow$}
\setbox42=\hbox{$\, \pi_2$}
\setbox4=\hbox{\kern -25mm{\box41 \raise 3mm\box42}}
\setbox5=\hbox{$G_{231}=(1-{y_1 \over x_1 })(1-{y_1 \over x_2 }) $}
\setbox6=\hbox{$(1-{y_1 \over x_1 })(1-{y_2 \over x_1 })= G_{312} $}
\setbox7=\hbox{$\pi_2 \, \downarrow$}
\setbox8=\hbox{$\downarrow \, \pi_1$}
\setbox9=\hbox{$G_{213}=(1-{y_1 \over x_1 })$}
\setbox10=\hbox{$(1-{y_1y_2 \over x_1x_2})= G_{132} $}
\setbox11=\hbox{$\kern 25 mm \pi_1 \, \searrow$}
\setbox12=\hbox{$\kern -25mm \swarrow \, \pi_2$}
\setbox13=\hbox{$G_{123 }\quad = \quad 1$}

$$ \matrix{   &\box2    \cr
\noalign{\vskip 1mm}
\box3 & &\box4 \cr
\noalign{\vskip 2mm}
\box5  & &\box6 \cr
\noalign{\vskip 3mm}
\box7  & &\box8 \cr
\noalign{\vskip 3mm}
\box9  & &\box10 \cr
\noalign{\vskip 3mm}
\box11 & &\box12 \cr
 &\box13 \cr  }$$

Let us also illustrate  the relations between the 
Yang-Baxter elements associated with the various Hecke algebras.

Take $\mu=(35142)$ as in example \ref{EX}. 

The coefficient of $T_2:=T_{13245}$ in
$Y^T_{35142}$ is the same as in the product
$$
\left( 1 + {x_3x_2^{-1} -1 \over q_1+q_2} T_2\right)
\left( 1 + {x_5x_1^{-1} -1 \over q_1+q_2} T_2\right)
\left( 1 + {x_5x_4^{-1} -1 \over q_1+q_2} T_4\right)
\left( 1 + {x_4x_2^{-1} -1 \over q_1+q_2} T_4\right)
$$
that is,
$$
{
(x_3x_5 - x_1x_2)(x_2x_4-q_1q_2(q_1+q_2)^{-2}(x_4-x_2)(x_5-x_4))
\over
(q_1+q_2 ) x_1x_2^2 x_4}
\ .
$$
The coefficient of $\bar\pi_2:=\bar\pi_{13245}$ in
$Y^{\bar\pi}_{35142}$ is the same as in the product
$$
(1+(1-{x_3\over x_2})\bar\pi_2)(1+(1-{x_5\over x_1})\bar\pi_2)
$$
that is,
$$
\left( 1- {x_3\over x_2} \right) +
\left( 1- {x_5\over x_1} \right)-
\left( 1- {x_3\over x_2} \right)
\left( 1- {x_5\over x_1} \right)
$$
$$
= 1 -{x_3x_5\over x_1x_2} \ .
$$
It is the specialization of the preceding coefficient
for $q_1=-1,q_2=0$.

The coefficient
of $\partial_2:=\partial_{13245}$ in $Y^\partial_{35142}$
is the same as in
$$
(1+(x_2-x_3)\partial_2)(1+(x_1-x_5)\partial_2) \ ,
$$
that is,
$$
x_1+x_2-x_3-x_5 \ .
$$
The coefficient of
$\sigma_2=(13245)$ in $Y^\sigma_{35142}$
is the same as in the product
$$
(1+(x_3-x_2)\sigma_2)(1+(x_5-x_1)\sigma_2)
(1+(x_5-x_4)\sigma_4)(1+(x_4-x_2)\sigma_4)
$$
that is,
$$
(x_1+x_2-x_3-x_5)(1+ (x_5-x_4)(x_4-x_2)) \ ,
$$
whose term of lowest degree corresponds to the preceding
case.


\end{document}